\documentclass{aupl}[final]

\usepackage{
amsfonts,
latexsym,
amssymb,
enumerate,
verbatim,
mathrsfs,
}

\newcommand{\labbel}[1]{\label{#1} [[{\bf #1}]]}  
\newcommand{\bibbitem}[1]{\bibitem{#1} [[{\bf #1}]]}  
\renewcommand{\labbel}{\label} \renewcommand{\bibbitem}{\bibitem}  

\newcommand{\sumbul}{\sum^\bullet}
\DeclareMathOperator{\nlim}{nlim}
\DeclareMathOperator{\mlim}{mlim}

\newtheorem{theorem}{Theorem}[section]

\newtheorem{proposition}[theorem]{Proposition} 
 
\newtheorem{corollary}[theorem]{Corollary}

\newtheorem*{claim*}{Claim}

\newtheorem*{theorem*}{Theorem}
\newtheorem*{proposition*}{Proposition}
\newtheorem*{corollary*}{Corollary}
\newtheorem*{lemma*}{Lemma}
\newtheorem*{scholion*}{Scholion}

\theoremstyle{definition}
\newtheorem{definition}[theorem]{Definition}

\theoremstyle{remark}
\newtheorem{remark}[theorem]{Remark}

\newtheorem*{remark*}{Remark}
\newtheorem*{remarks*}{Remarks}

\newtheorem{observation}[theorem]{Observation}

\newtheorem*{observation*}{Observation}

 \allowdisplaybreaks[1]

\begin{document}

\title{Series of combinatorial games}
 
\author{Paolo Lipparini} 
\urladdr{http://www.mat.uniroma2.it/~ lipparin}
\address{Dipartimento  di Matematica\\Viale della  Ricerca
 Combinatoria\\Universit\`a di Roma ``Tor Vergata'' 
\\I-00133 ROME ITALY (currently retired)}


\thanks{Work performed under the auspices of G.N.S.A.G.A. 
The author acknowledges the MIUR Department Project awarded to the
Department of Mathematics, University of Rome Tor Vergata, CUP
E83C18000100006.}

\begin{abstract}
We present a  definition 
for the sum of a sequence of combinatorial games.
This sum  coincides with the classical sum in the case of a converging
sequence of real numbers and with the infinitary natural sum
in the case of a sequence of ordinal numbers.

We briefly discuss other possibilities, such as the string limit,
some ``magical'' variants of Hackenbush,  as well as
``Dadaist'' infinite sums, which allow transfinite runs,
while still being loopfree.   
\end{abstract} 

\maketitle

\section{Introduction} \labbel{intro}

The problem of evaluating infinite sums has a long history,
in a sense---in hindsight---starting as early as  
Zeno's Achilles and the tortoise paradox \cite{stan}. See also \cite{har}.

Conway's surreal numbers \cite{onag} generalize at the same time
real numbers and ordinals. Surreal numbers form a Field,  have 
been studied from different points of views \cite{Al,A,onag,E,G,S}
and have been recently found useful in making progress
to the solution of  an old problem by Skolem \cite{BM}. 

Needless to say,  limits and series play a fundamental
role in analysis. Limits and infinite sums are also present 
in the theory of ordinals \cite{bach}. 
See \cite{L,W}  for further references.
Various  notions of ``surreal limits'' or 
infinite sums appear, among others, in \cite[p. 40]{onag}, 
\cite{RS},  \cite[Definition 3.19]{S} and \cite{M,LM}.
In passing, notice that the sum defined in \cite[p. 40]{onag}
does not necessarily coincide with the ordinal sum, for example,
$ \omega +  \omega^2 + \omega ^3 + \dots $ is not defined in the sense
of \cite[p. 40]{onag}. 

 Conway \cite{onag}
also introduced a general notion of ``games''
of which ``numbers'', i.e., surreal numbers,  
are a proper subclass. 
In passing, it is a curious fact that positive infinitesimal
games, actually, small games (which however are not numbers)
occur already at Day 2 \cite[Section II.4 and p. 404]{S}.

To the best of our knowledge, no
notion of infinite sum or limit has ever been proposed for arbitrary 
combinatorial games.
Here we introduce a notion of  series 
for arbitrary combinatorial games and evaluate it in
a few cases. In particular, the sum of a (classically)
converging series of reals 
has the same value also in a game theoretical sense.
In the case of ordinal numbers, the sum we introduce
coincides with the so-called \emph{infinitary natural sum}.
See \cite{L} for details and further references on ordinal sums. 

In conclusion, what we propose is a purely game-theoretical definition 
of an infinite sum which coincides with classical and useful notions
in at least two significant cases.

\section{Series of games} \labbel{sgam}

We assume familiarity with the basic notions of 
combinatorial game theory.
A \emph{combinatorial game}, or simply \emph{game} for short,
  is a position in a possibly infinite
two-player game with perfect information, no chance element, no draw
and no infinite run, as studied, for example, in \cite[Chapter VIII]{S}.
Games are possibly \emph{partizan}
and we always assume \emph{normal} (not mis\'ere) \emph{play}.  
The two players are called \emph{Left}
and \emph{Right} and in each play on some game
they can move either as the first or 
the second player. We sometimes capitalize
\emph{First} and \emph{Second} in order
to denote the player who plays first or second
in some specified play, no matter whether she  is Left or Right.

\begin{definition} \labbel{sum} 
Given a sequence $(G_i) _{i  \in \mathbb N} $
of games, we are going to define a game 
$\sum _{i  \in \mathbb N} G_{i}$, sometimes written
as $G_0+G_1+G_2+ \dots $ or, possibly,
  $G_0+G_1+ \dots + G_i+ \dots $, etc.
 The game $\sum _{i  \in \mathbb N} G_{i}$
will be called an \emph{infinite sum} 
and sometimes a \emph{series}, in order to avoid confusion
when comparing it with some finite sum, e.~g.,
some partial sum. Thus we do not reserve the word
``series'' only to the classical notion of analysis 
(we will see in Section \ref{realsec} that the notion
we are going to define actually incorporates the classical
notion).   

Rather informally, a play on 
$\sum _{i  \in \mathbb N} G_{i}$ goes as follows.
If First has no move on any $G_i$, First looses the game.
Otherwise, First chooses some natural number $n$
and she makes a move $G^F_i$ on some $G_i$  with $i \leq n$. 
The resulting position on $\sum _{i  \in \mathbb N} G_{i}$  
after  First move is 
 $G_0+ \dots + G^F_i+ \dots + G_n + G _{n+1} + \dots  $.

Then Second 
 chooses some natural number $m \geq n$ and   he
makes a move $G_j^S$ on some $G_j$  with $j\leq m$
(of course, moving on $G^F_i$, not on $G_i$, if $j=i$).
Then the play continues on the finite sum
 $G_0+ \dots + G^F_i+ \dots + G_n + \dots +  G _{j}^S + \dots + G_m $
(possibly with the double subposition addendum $G^{FS}_i$, instead).
   
For short, each player must 
 choose some index as soon as she plays, and then she 
has to move on some game with the same or  a smaller index.
After both players have made their choices
of the indexes, the play continues on the
finite sum of the  games
up to the larger chosen index.
The same description applies to \emph{runs},
ie, when the same player is supposed to make
 two or more consecutive moves
(a \emph{play} is an alternating run;
the development of Conway theory
requires the general definition of runs).
Notice that, after her choice of $n$,
and until Second moves, First
can play only on $G_0+ \dots +G_n$, which is  a finite sum of games,
hence, for every run, either the run ends after a finite number
of moves, or, sooner or later, Second moves. After this, the game 
definitely proceeds 
on a  finite sum with a fixed index $m$, hence 
the game eventually terminates.
Thus 
$\sum _{i  \in \mathbb N} G_{i}$
has no infinite run.

We now define 
$\sum _{i  \in \mathbb N} G_{i}$  formally.
We need to introduce auxiliary games
such as  
$ \left( \sum _{i  \in \mathbb N} H_{i}\right ) / {L,n}$,
with the intended meaning that ``Left has
already made the choice of her index 
and she has chosen $n$''. 
For clarity, we may possibly write
 $(H_0+H_1+H_2+ \dots ) / {L,n}$ in place of
$ \left( \sum _{i  \in \mathbb N} H_{i}\right ) / {L,n}$.
Thus
\begin{align*}  
&\big( \sum _{i  \in \mathbb N} H_{i}\big) / {R,n} 
=
\\
&\big\{ \, 
H_0+ \dots +H^L_j+ \dots + H_m  \  \big| \ 
  (H_0+ \dots +H^R_i+ \dots + H_n + \dots ) / {R,n}  \, \big\}, \text{ and}
\\
&\big( \sum _{i  \in \mathbb N} H_{i}\big) / {L,n} 
=
\\
&\big\{ \, (H_0+ \dots +H^L_i+ \dots + H_n + \dots ) / {L,n} \  \big| \ 
H_0+ \dots +H^R_j+ \dots + H_m \ \big\}
\end{align*} 
where $m$ is always intended to be $\geq n$.  
For each fixed $n$, the above definition 
is by transfinite induction on the natural sum of the 
birthdays of $H_0, \dots, H_n$.
Indeed, after some move on $H_i$,  
 the natural sum of the 
birthdays of, say,   $H_0, \dots, H^L_i, \dots,  H_n$
is strictly smaller than the natural sum of  the 
birthdays of $H_0, \dots, H_i, \dots,  H_n$,
hence the induction carries over.
Notice that 
$H_0+  \dots + H_m $ 
might in principle have a birthday as large as any ordinal,
but $H_0+  \dots + H_m $ is well defined in any case and the 
induction does not involve its birthday.
Finally,
\begin{align*}
 \sum _{i  \in \mathbb N} G_{i} =
& \{ \, (G_0+ \dots +G^L_i+ \dots + G_n + \dots ) / {L,n}\  \big | \ 
 \, 
\\ & 
(G_0+ \dots +G^R_i+ \dots + G_n + \dots ) / {R,n} \} 
  \end{align*}     
\end{definition}   

The operation introduced in Definition \ref{sum}
is not invariant under Conway equivalence of the summands.
See Remark \ref{serve} in the appendix.  

We now state some simple facts about
the above notions. As usual in the classical theory
of series, if $\sum _{i  \in \mathbb N} G_{i}$
is an infinite sum, the corresponding \emph{partial sums} 
are the finite sums of the form
$G_0 + \dots + G_{h-1}$, for some $h \geq 0$
(the notation allows the possibility for an empty sum,
the $0$ game, 
however, here we will have no use for it).   
 
Recall that $G \geq 0$ if Left has a winning strategy
when playing second on $G$, and 
$G > 0$ if Left has a winning strategy
whoever starts the game.
Recall that a  \emph{number} 
is a game $G$ such that $H^L < H^R$,
for every (possibly improper) subposition $H$ of $G$. 
A game is \emph{dicotic} if both players have a move
on every nonempty subposition of G. 

\begin{observation} \labbel{obs}   
Suppose that $(G_i) _{i  \in \mathbb N} $
is a sequence of games.
  \begin{enumerate}    
\item
If each $G_i$ is impartial (a number, dicotic, $>0$, $\geq 0$,  $<0$, $\leq 0$),
then  $\sum _{i  \in \mathbb N} G_{i}$ is
impartial (a number, dicotic, $>0$, $\geq 0$,  $<0$, $\leq 0$).

\item
If there are infinitely many partial sums which are
$\geq 0$ (resp., $\leq 0$), then Left (resp., Right)
wins when playing second on  $\sum _{i  \in \mathbb N} G_{i}$.
 In particular, if both eventualities hold, then 
$\sum _{i  \in \mathbb N} G_{i}$ is second winner.

\item
If each $G_i$ is an ordinal, then
 $\sum _{i  \in \mathbb N} G_{i}$ is the supremum
(in the ordinal sense)
of the partial sums (which, in the ordinal sense, are
intended as \emph{natural sums} \cite{bach}).
 
\item
If the $G_i$ are eventually Conway zero games,
ie, second winner games,
then $\sum _{i  \in \mathbb N} G_{i}$ is
Conway equivalent to the finite sum of the nonzero summands.

\item
 An infinite sum of 
  impartial games is a 
second winner game, if infinitely many summands
are first winner (notice that if only finitely many summands
are first winner, the outcome is given by (4),
since an impartial game is either first  or second winner).
 \end{enumerate}
\end{observation} 

\begin{proof}
(1) We prove the statement in the case of numbers, all the other statements 
are much simpler. We first show that 
$\big (\sum _{i  \in \mathbb N} H_{i}\big) / {R,n}$
is a number.
We have to show that 
Left always wins playing on
$  (H_0+ \dots +H^R_i+ \dots + H_n + \dots ) / {R,n} 
- H_0-H_1 \dots -H^L_j \dots - H_m  $.
Since $H_i$ and $H_j$ are numbers by assumption,
$H^R_i > H_i$ and
$H^L_j < H_j$, that is, 
$-H^L_j > - H_j$
so Left wins moving first
simply by choosing $m=n$ in   
 $(H_0+ \dots +H^R_i+ \dots + H_n + \dots ) / {R,n}$.
Again since all the summands are numbers,
if Right plays first, he increases the value of the summand
he moves on; a fortiori,  Left wins moving second.
Similarly, $\big ( \sum _{i  \in \mathbb N} H_{i}\big) / {L,n}$
is a number.

We also need to show that 
$(G_0+ \dots +G^L_i+ \dots + G_n + \dots ) / {L,n}
< (G_0+ \dots +G^R_i+ \dots + G_{n'} + \dots ) / {R,n'}$,
that is,
$(G_0+ \dots +G^R_i+ \dots + G_{n'} + \dots ) / {R,n'}
+ (-G_0-G_1 \dots -G^L_i+ \dots - G_n + \dots ) / {R,n} 
>0$.
This is similar to the above arguments; 
here Left has to choose her $m$'s both in the first and in the second
sum; it is enough that she chooses the same $m$
in both games, as soon as she plays there, and that $m \geq n, n'$. 

(2) After First has chosen $n$, let Second choose some 
$m \geq n$ such that 
  $G_0+  \dots  + G_{m} \geq 0$, respectively,
$\leq 0$. 

(3) The surreal sum of two ordinals is their natural sum.
Right cannot move on an ordinal, and 
a Left move on $\sum _{i  \in \mathbb N} G_{i}$
simply amounts to choosing some $m$ and moving
on the corresponding partial sum. This is exactly the same
as moving on the ordinal limit of the partial sums. 

(4) Let $G _{i_1}, G _{i_2}, \dots \  $
enumerate the nonzero games.
 Second has a winning strategy on
$(\sum _{i  \in \mathbb N} G_{i})- G _{i_1} - G _{i_2}- \dots  $
by choosing a sufficiently large index as soon as 
he plays on the infinite sum. He uses the mirror-image strategy
for the nonzero games and the winning strategy as a second player
on the 
(possibly empty set of)
Conway zero games.

(5) Let First choose some $n$, moving to
$  (G_0+ \dots +G^F_i+ \dots + G_n + \dots ) / {F,n} $.
If $G_0+ \dots +G^F_i+ \dots + G_n$ is first winner,
then Second wins by choosing $m=n$.
Otherwise,   
$G_0+ \dots +G^F_i+ \dots + G_n$ is second winner,
and then Second wins by choosing $m > n$ in such a way that
the ``new'' set of games contains exactly one first winner game;
then he moves on that game adopting the winning strategy
for the first player.  
\end{proof}

The infinite ordinal sum described in (3)
has been extensively studied in \cite{L}. See
\cite[Definition 2.3]{L}.

\section{Series of real numbers} \labbel{realsec} 

Recall that each real number can be considered
as a game; we will 
always consider a real number as a surreal number 
in
\emph{surreal canonical representation} 
\cite[Theorem 2.8]{G}
(not to be confused with the \emph{canonical form}
of a short game \cite[Section II.2]{S}).
The surreal canonical representation 
of a surreal corresponds
to its sign expansion, the sign sequence
which gives an equivalent game expressed as
a single, possible transfinite, stack of Red Blue
Hackenbush.

\begin{theorem} \labbel{reals}
Suppose that  $\sum _{i  \in \mathbb N} r_{i}$
is a series of real numbers, converging
(in the classical sense of analysis) to 
some real number $r$. Then
 $\sum _{i  \in \mathbb N} r_{i} $
is Conway equivalent to $  r$.
 \end{theorem}

  \begin{proof}
We have to show that 
 $-r + \sum _{i  \in \mathbb N} r_{i} $
is a second winner game.
As we mentioned, we assume
that all the numbers under consideration
are expressed in surreal canonical representation. 
By Observation \ref{obs}(4) we may also  assume that
$( r_i) _{i \in \mathbb N} $ is not eventually 
$0$. 

Suppose that 
Right moves first (the other case
is treated symmetrically).

(a) If Right moves
on $\sum _{i  \in \mathbb N} r_{i}$,
he chooses some $n \in \mathbb N$ and moves on some
 $r_i$,  getting a subposition $r_i^R > r_i$,
since $r_i$ is a (surreal) number.
Thus the series $r_0+ \dots + r_i^R+ \dots  + r_n + \dots $
classically converges to  $s= r+r_i^R - r_i$ and  $s > r$.
By convergence of the latter series,  there is some $m \geq n$ such that 
 $H=r_0+ \dots + r_i^R+\dots  +  r_n + \dots + r_m >r$.
That is, $-r+ H =- r + r_0+ \dots + r_i^R+\dots  +  r_n + \dots + r_m > 0$
and this means that Left has a winning strategy on $-r + H$. 
If the winning strategy requires Left to move on $-r$,
let her  move there. Then Right can move 
either on $(-r)^L$ or on  $r_0+ \dots + r_i^R+\dots  +  r_n$,
but in the latter case he cannot change the value of $n$.
As soon as Left strategy involves moving on $H$,
let her  choose the natural number $m$.  
 Then the game continues on 
(the remaining options of)
 $-r + H$, so Left wins.

(b) Otherwise, Right moves on $-r$ and
we have $(-r)^R > -r$, again by a property
of numbers. 
Let $-s= (-r)^R$,
thus $s < r$. 
Moreover, 
$(-r)^R$  has finite birthday
since a real number in surreal canonical representation
has birthday $\leq \omega$
and  options have always strictly smaller birthday. 
So there is some $h \in \mathbb N$ 
such that Right can perform at most $h$
consecutive moves on  $(-r)^R$.

Since $\sum _{i  \in \mathbb N} r_{i} = r$
in the classical sense and $s < r$,
where $r$ and $s$ are real,   there is $n' \in \mathbb N$ 
such that $r_0 + r_1+ \dots + r_m > \frac{r+s}{2} $,
for every $m \geq n'$. 
Since the series is converging, $\lim _{ n \in \mathbb N} r_n =0 $,
thus there are $n'', p \in \mathbb N$    
such that  $|r_{\ell} |< \frac{1}{2^p} < \frac{r-s}{2(h+1)}$,
for every $\ell \geq n''$.

Let Left move on $\sum _{i  \in \mathbb N} r_{i}$
and choose an $n$ such that $n \geq n'$,
$n \geq n''$
and such that, moreover, there are at least $h+1$  distinct 
indexes $ \ell_0, \ell_1, \dots, \ell_h$
such that the games
$ r _{\ell_0}, r _{\ell_1}, \dots, r _{\ell_h}$ are nonzero
and moreover   $n'' \leq \ell_0\leq n $, \dots, $ n'' \leq \ell_h \leq n$.
As we mentioned at the beginning of the proof,
we may assume that there is such a sequence
of nonzero games in view of
Observation \ref{obs}(4). 
We claim that Left has a winning
strategy choosing $r _{\ell_0}$ for her move.
If $r _{\ell_0} >0$, she can move to $0$ on $r _{\ell_0} $,
since games are in  surreal canonical representation,
so the value of the total game decreases by $r _{\ell_0}$. 
If $r _{\ell_0} <0$, she can 
decrease the value of the game by 
less than $\frac{1}{2^p}$, since
 $|r _{\ell_0} |< \frac{1}{2^p}$, so the 
sign sequence of 
$r _{\ell_0} $ starts with a minus sign followed by
at least $p+1$ plus signs, so she can remove 
the $p{+}1$th plus sign and whatever follows. 

(b1) If Right's second move is on the infinite sum,
he must choose some $m \geq n$.
Since  $n \geq n'$, we have
$r_0 + r_1+ \dots + r_m > \frac{r+s}{2} $ for the original games.
However, Left has played on $r _{\ell_0} $,
so we need to compute $r_0 + \dots r _{\ell_0}^L+ \dots + r_m
 > \frac{r+s}{2}- \frac{1}{2^p} >
 \frac{r+s}{2} - \frac{r-s}{2(h+1)} \geq
\frac{r+s}{2} - \frac{r-s}{2} = s $.
Since $-s=(-r)^R$, Left wins on
$(-r)^R + r_0 + \dots r _{\ell_0}^L+ \dots + r_m$.  

(b2) Otherwise, Right second move is again on 
$(-r)^R $ and $(-r)^{RR} > (-r)^R $.
Left can go on moving on 
$r _{\ell_1}$,  as above, 
contributing a negative difference
of less than $\frac{1}{2^p}$ 
to the value of the total game.
Since Right can move only 
$h$ times on $(-r)^R $,  
Left has a sufficient number of games 
$r _{\ell_k}$ to move on.

Eventually, Right is forced to move on the infinite sum
and he must choose some $m$.
If $K$ is the sum of the   
options remained so far, then, arguing as above, 
$K > \frac{r+s}{2}- \frac{h+1}{2^p} >
 \frac{r+s}{2} - \frac{(r-s)(h+1)}{2(h+1)} =
\frac{r+s}{2} - \frac{r-s}{2} = s $.
Hence Left wins on
$(-r)^{R} + K $,
hence also on 
$(-r)^{RRRR\dots} + K $,
since any move by Right 
on $(-r)^{R}$
increases the value of the game.
 \end{proof}

\section{Further remarks} \labbel{fur}

\begin{remark} \labbel{nolim}    
(a) As we have seen in the previous sections, the infinite sum
introduced in Definition \ref{sum} has some good 
properties. However, it misses an important property 
a notion of limit is usually supposed to share, namely,
it is not always true that the outcome of a sum  depends only on the 
``tail'' of the sequence of the partial sums. In detail, it could be expected
that if $\sum _{i  \in \mathbb N} G_{i}$
and $\sum _{i  \in \mathbb N} H_{i}$
are two series and
 there is some $h$ such that, for every $k \geq h$,  the two partial sums
with  $k$ elements   have the same value,
then the series have the same value.

This is not always true for 
the notion introduced in Definition \ref{sum}.
For example,    
 $-1+1+1+1+ 1+\dots = \omega -1 $ and 
$0+0+1 + 1+1+\dots  = \omega $,
though, from $h=2$ onward, the partial sums have the same value.  

We do not know whether 
there is some variation on Definition \ref{sum}
such that the above desirable property is satisfied. 
Of course, the above ``tail property'' holds in many significant cases,
e.~g., sums of ordinals, or converging sums of real numbers
in canonical representation,
by Observation \ref{obs}(3) and Theorem \ref{reals}.   

(b) Remark (a) suggests that a related definition is quite natural.
Define an infinite sum as in Definition \ref{sum},
except that each player, once and for all,
can choose subsets $I$ and $J$  of $\mathbb N$, rather
than simply natural numbers $n,m$.
 After both choices have been made, the
game is played on the finite sum 
$\sum _{i  \in I \cup J} G_{i}$, while similar rules
as in Definition \ref{sum} are applied
just after the first choice. 

With the resulting infinitary operation, for real summands,
the notion of \emph{absolute convergence}  is recovered.
There is no essential difference with respect to the proof
of Theorem \ref{reals}.
Moreover, this modified operation is 
invariant under permutations of $\mathbb N$.  
When restricted to ordinals, the two kinds of
sums coincide \cite{L}.

(c) Series are connected with limits, so we can 
associate to Definition \ref{sum} a notion of limit
for sequences of combinatorial games.
In detail, if $( H_i) _{i \in  \mathbb N} $ 
is a sequence of games, 
we may define
$\lim _{i \in \mathbb N}  H_i = \sum _{i \in \mathbb N} (H_i - H_{i-1})$, 
where $H_{-1} = 0$.

This translates to a direct definition.
For short, First makes some move $H_i^F$  on some 
 $H_i$, then Second might continue moving on 
$H_i^F$ or  chose some $j >i$. 
If Second chooses the latter option, she moves on 
  $H_j -  H_i + H_i^F$ and the play continues there.

This limit has some natural properties,
corresponding to Observation \ref{obs}.
If the sequence  $( H_i) _{i \in  \mathbb N} $
is constant  from some point on, the limit is
Conway equivalent to
this constant value.
If  $H_i \geq 0$ for infinitely many indexes $i$,
then the limit is $\geq 0$. 

As in (a) above, this limit generally depends
also on the first items of the sequence, and not only on
its ``tail'', though there are significant cases in which
this dependence does not occur.

(d)
The sum introduced in Definition \ref{sum} 
satisfies really few associativity properties.
For example,  $-1+1+1+1+\dots= \omega-1$ but
  $(-1+1)+1+1+\dots = \omega$.

Of course, under really modest assumptions,
there are many forms of associativity
which  cannot be expected to hold
for an \emph{everywhere  defined} infinitary operation
 \cite[Example 2.5]{L}. See \cite[Chapter 5]{W} 
for an ample discussion of incompatibilities of desiderata
for infinite sums.

(e) The behavior of
$\sum $ is not always good with respect to addition, either.
 $1+1+1+1+\dots= \omega$ but also
 $2+2+2+2+\dots= \omega$.
Hence $\sum _{i  \in \mathbb N} G_{i} +
 \sum _{i  \in \mathbb N} H_{i}$
is not necessarily equal to 
$\sum _{i  \in \mathbb N} (G_{i}+H_i)$.

The same example shows that it is not always the case
that, for surreals, $ \lambda \sum _{i  \in \mathbb N} G_{i}
=\sum _{i  \in \mathbb N} \lambda G_{i}$.
Recall that, for arbitrary games, the product is not even  invariant
under Conway equivalence. 

(f) Definition \ref{sum} can be extended with no 
particular difficulty to sums indexed by an ordinal. 

Define $\sum _{ \beta  < \alpha } G_{ \beta }$
by transfinite induction on $\alpha$ as follows.
We set 
\begin{equation*}
\sum _{ \beta  < 0 } G_{ \beta } = 0
\end{equation*}    
and
\begin{equation*}  
\sum _{ \beta  < \alpha+1 } G_{ \beta }
= ( \sum _{ \beta  < \alpha } G_{ \beta }) + G_ \alpha .
\end{equation*}     

For $\alpha$ a limit ordinal, a play on 
$\sum _{ \beta  < \alpha } G_{ \beta }$ goes as follows.
First chooses some $\beta < \alpha $ and moves on
some $G_ \gamma $ with $\gamma < \beta $.
At Second turn, she chooses some $\delta$ such that  
$\beta \leq \delta < \alpha $ and 
 the play continues on (the remaining
options of) $\sum _{ \beta  < \delta  } G_{ \beta }$,
which has been already inductively defined, since
$\delta < \alpha$.

When $ (G_{ \beta })_{ \beta  < \alpha }$ 
is a sequence of ordinals, $\sum _{ \beta  < \alpha } G_{ \beta }$
turns out to be the same as the sum studied in 
\cite{trans}. See Definition 3.1 therein and see the quoted source
for historical remarks and further references.   

We expect that if the surreal $s$ is expressed in Conway normal form
as  $\sum _{ \beta  < \alpha } \omega ^{a_ \beta } r_ \beta  $
then $s$ is also the value of the sum computed according
to the above rules, but we have not
fully checked every detail, yet. 
 \end{remark}

\section{Other possibilities} \labbel{other}

\subsection{The string limit} \labbel{string} 

One can introduce the notion of the  limit of
a sequence of strings of symbols.
Here we consider the simplest case of an $\mathbb N$-indexed 
sequence of strings, the more general case of ordinal
indexed sequences of strings presents no essential
difference. See \cite[Section 4]{LM} for even more
general definitions.

A \emph{(transfinite)  string} is an ordinal-indexed
sequence of elements from an arbitrary
set $S$, considered as a set of symbols.
For example, when considered as a \emph{sign-sequence},
a surreal number is just a string 
of elements from $S= \{ +. - \} $.
   
If each $s_i$ is a string (and they are possibly indexed by
distinct ordinal numbers) 
the \emph{string  limit} of the sequence $( s_i) _{i \in \mathbb N} $ 
is the longest ordinal-indexed sequence $s$ 
such that
  \begin{enumerate}[(I)] 
  \item
If $s$ has successor length
(including length $0$, that is, $s$ is the empty string),
then  there is $n \in \mathbb N$ such that 
$s$ is the initial segment of every 
$s_m$, for $m \geq n$.      
\item
If $s$ has limit length,
then, 
 for every proper initial segment $s'$  of $s$,
there is $n \in \mathbb N$ such that 
$s'$ is the initial segment of every 
$s_m$, for $m \geq n$.      
   \end{enumerate} 

The definition is justified by the following observations.
If two strings of lengths $\alpha$ and $\beta$ 
satisfy one of the above conditions, they coincide
up to $\min \{  \alpha , \beta \} $.
If there is no maximal string satisfying (I),
then any string satisfying   (I) can be extended to some longer
string still satisfying (I). The common
prolongment  of all such strings
satisfies (II).

For example, the string limit of 
$+$, $+-$, $+-+$, $+-+-$ \dots 
\ is the infinite string 
$+-+-+- \dots$\ (here we are denoting strings
with juxtaposition).
The string limit is not necessarily as long
 as the inferior limit of the lengths of the strings in 
the sequence.   The limit of 
$+$, $-$, $+-$, $-+$,  $+-+$, $-+-$ \dots 
\ is the empty string.
Intermediate cases might occur,
for example,
the string limit of 
$+$, $+-$, $+++$, $+-$,  $+-++$, $+++++$,
$+-$, $+-++++$ \dots \    is
$+$, since the first place remains fixed, while the second 
place alternates.    

As a final example, the string limit of the sequence 

$+----\dots\  -$, 

$++---\dots\  -$, 

$+++--\dots \ -$, 

$++++-\dots \ -$, 

is 

s=$+++++\dots $,

an $\mathbb N$-indexed sequence of $+$. 
Indeed, any proper initial segment of $s$,
a sequence of $n$ pluses, is the initial segment
of all the strings $s_m$, for $m \geq n$.    
Each string $s_i$ contains a minus sign at the 
$ \omega$th place,  but this minus sign does not
belong to $s$, since    
$+++++\dots -$
is never an initial segment of some $s_i$. 

The string limit of strings from $S= \{ +. - \} $,
that is, surreal numbers expressed as sign sequences,
has been extensively studied in \cite{LM}. 
There it has been proved that the string limit
of a converging sequence of real number
coincides with its limit, modulo an infinitesimal.
Moreover, the string limit of a non decreasing
sequence of surreal numbers is $\geq$ than
each element of the sequence.\footnote{In passing, note that the statement
on paragraph -2 on p. 9 in \cite[v. 2]{LM} about some purported ``continuity''
of the string limit does not always hold.}
 
The limit from \cite{LM} can  be used also in order to 
provide a simple proof \cite{clean} of the  Fundamental  Existence  Theorem
for separating elements of sets of surreals \cite[Theorem 2.1]{G}.
The simple proof from \cite{clean} has been formalized and
verified by Metamath with contributions by Scott Fenton \cite{Met}.

The above notion of  string limit is probably interesting
also in settings distinct from game theory.
In any case, it can be applied to 
any game which can be represented as a single 
(possibly transfinite) stack of Green-Red-Blue Hackenbush.
For example, the string limit of $*n$
(an $n$-stack of green edges) is $* \omega $,
an $\mathbb N$-indexed stack of green edges.  
 
Notice that, contrary to the limits
considered so far, the string limit satisfies 
the ``tail property'': modifying
finitely many elements of a sequence
does not change the limit of the sequence.
Cf. Remark \ref{nolim}(a).

\subsubsection{Some ``magic'' variations of
Hackenbush} \labbel{hack} 

In Hackenbush one could also introduce some further kinds of 
``magic'' edges, say, a \emph{dark green}  (\emph{dark blue, dark red})
 edge which can be removed only
if it has  one free 
vertex; so that 
if it lies in a single stack, it can
be removed only if it lies at the top of the stack.
Of course, such a magic edge  ``flies away'' 
if it becomes disconnected to the ground---in a stack, if some edge below 
 is removed. 

Another possibility is a 
\emph{greenish blue edge} 
which can be always removed by Left,
but can also be removed by Right 
if it has  one free 
vertex. Thus $ \uparrow $
can be realized by a  
red edge atop such a 
greenish blue edge.

Still another very magical possibility,
consider a
\emph{super-red edge}  which, when removed, turns into
green all the edges 
which were connecting it to the ground.
Again,
$ \uparrow $
can be realized: by such a  
super-red edge atop  a 
blue edge.

Any sequence of stacks containing
ordinary and magic ``edges''
has a string limit  defined as above,
though we expect that this limit may provide rather strange 
results.
 There are obvious possible variations
in the definition of the string limit. For example, 
when the strings in a sequence contain
alternately a  blue edge 
and a greenish blue edge at the same position,
we might allow the limit to be defined at that position
and to be a greenish blue edge.

One can also generalize the string limit 
by taking limits modulo filters
(to be defined at some place, instead of ``eventually constant''
at that place,
it is required that the string is constant at that place modulo $F$---thus
we get the string limit exactly when $F$ is the Fr{\'e}chet 
filter of cofinite subsets of $\mathbb N$),
and also by considering a topology on the set $S$ of
symbols (at each place, we consider topological convergence,
rather than being eventually constant---note that the latter
is convergence modulo
the discrete topology).
One might also simultaneously combine the above generalizations.

Full details appear in Chapter 4 of D.
Battistelli thesis \cite{Ba}.

\subsection{A natural limit} \labbel{natur} 

There is a quite natural notion of a limit of games.
If $(G_i) _{i  \in \mathbb N} $ is a sequence of games,
the \emph{natural limit} 
 $\nlim  _{i  \in \mathbb N} G_i $
is defined as follows.
First chooses some $n \in \mathbb N$ 
and some option $G^F$  which is present in 
\emph{each} $G_i$, for $i \geq n$.
If no such option exists, First loses.
Otherwise, the play continues on 
$G^F$.

If the $G_i$ are surreal numbers expressed 
in canonical representation as sign
sequences, then nlim coincides with the
string limit. Similarly for a sequence of games
which are expressed as a single Hackenbush stack.
If the sequence of the $G_i$ is eventually constant, then
 the natural limit coincides with this constant game.
Moreover, contrary to many other limits considered here,
the ``tail property'' is satisfied: 
the natural limit does not depend on any finite initial part
of the sequence.

On the other hand, the natural limit is not invariant under Conway
equivalence and heavily depends on the form of the games
under consideration.
If we consider natural numbers as sign sequences,
namely, $n+1= \{ \, 0, 1, \dots, n \mid  \,\} $,
then $\nlim  _{n \in \mathbb N} n= \omega $.  
If, instead,  we represent natural numbers in
canonical form as
 $n+1= \{ \, n \mid  \,\} $,
then no player can move on the limit,
hence $\nlim  _{n \in \mathbb N} n= 0 $.

\emph{A variant.} The above paradoxical aspects
of the natural limit can be mitigated by a modified definition.

If $(G_i) _{i  \in \mathbb N} $ is a sequence of games,
the \emph{monotone limit} 
 $\mlim  _{i  \in \mathbb N} G_i $
is defined as follows.

If Left plays first, she chooses some $n \in \mathbb N$ 
and some option $G_n^L$ such that,
for every $i \geq n$, there is some option
$G_i^L\geq G_n^L$.
If no such a sequence of options exists, Left loses.
Otherwise, the play continues on 
$G_n^L$.
 
Symmetrically, if Right plays first, he chooses some $n \in \mathbb N$ 
and some option $G_n^R$ such that,
for every $i \geq n$, there is some option
$G_i^R\leq G_n^R$.
Then the play continues on $G_n^R$.

So if, again,  we represent natural numbers in
canonical form as
 $n+1= \{ \, n \mid  \,\} $,
then  $\mlim  _{n \in \mathbb N} n= \omega  $.  

As usual, as soon as we can introduce some notion of limit,
we also have the notion of an infinite sum,
just taking the limit of the partial sums.
Thus the string limit induces an infinitary sum
 operation on the surreal numbers, as already noticed
in \cite[Problem 5.4(a)]{trans}, and all the  limits
considered in the present subsection 
induce an infinitary sum
 operation on the class of all games.

\subsection{Further possibilities} \labbel{furfur}

The rule for \emph{conjunctive} sums  \cite[p. 41]{S}
can be obviously extended to any 
infinite set of games (any player must move on \emph{every}
summand).
The rule for \emph{ordinal} sums  \cite[p. 41]{S}
can be obviously extended to any ordinal-indexed sequence
of games
(moving on a summand
annihilates all the following summands). 

The \emph{sequential}  rule can also be applied, but using
the ordinal rule at limit stages.
For example, in 
$G _{ \omega }\rightarrow \dots \rightarrow G_n \rightarrow \dots \rightarrow 
 G_1 \rightarrow G_0   $ any player
should move on  
$G _{ \omega } $, unless $G _{ \omega } $ is $0$.   
If $G _{ \omega } $ is $0$, the player chooses some $n$
and moves on $G_n$. 
The play continues on $G_n$
until it is $0$; then the play continues on 
$G _{n-1} $ and so on.

Also the rules for the \emph{side} sum of two games 
\cite[p. 41]{S} can be
generalized, say, to a $\mathbb Z$-indexed sequence of games.
A Left move on some game annihilates
all the games on the left (that is, all the games
with a larger index) 
and a Right move on some game annihilates
all the games on the right.

In all the above cases the sum
 produces a combinatorial game 
(ie, without loops and without infinite runs).

If instead we use the \emph{selective}
rule (a player
can move in any number of components, but at least one),
 we get infinite runs, a situation similar to the one we will
describe in the next subsection for the classical \emph{disjunctive} 
rule (each player moves on exactly one summand).
 With the \emph{continued conjunctive} rule 
(each player must move in all non $0$  summands; the play
 ends when all summands are $0$)
we still get infinite runs, but in this case the maximal length
of a run is $ \omega$.
The main problem for all kinds of sums described in the present
paragraph, as well as in the following subsection, 
is to define the notion of an alternating run, ie,
to establish which player is supposed to move at some
limit step. 

Still another possible game-theoretical definition 
of the sum of a sequence of games---so far, limited
to ordinals---will be presented in \cite{anoth}.

\subsection{``Dadaist'' infinite games} \labbel{dada} 

There is another possible definition for a sum
$\sum _{i \in I} G_{i}$ 
of combinatorial games, where $I$ 
is an arbitrary, possibly infinite,
 set of indexes. Simply let each player, at her turn,
move on some game $G_i$, or
in the subposition which has remained. This rule gives rise---of 
course---to a game with possibly infinite 
(actually, transfinite) runs, hence not a combinatorial game
in the strict sense. After, say, $ \omega$ moves and for each
fixed $i$, only a finite number of moves have been played on 
the combinatorial game $G_i$,
hence some played subposition of $G_i$, possibly
$0$,  remains fixed from some point on.   
If at the $ \omega$th step we take all such remaining subpositions,
 we still have a possibly infinite game on which
the play can go on (unless all the summands have been
already exhausted).  We can do the same at each
subsequent limit step. If $I$ is a set, then, 
 assuming the axiom of choice \cite{J}, each play eventually
terminates after a run of length $<|I |^+$ 
(the successor cardinal of $| I |$, thought of as an ordinal),
 since we have $<|I |^+$  subpositions in total. 
Thus, while some proper subposition of $\sum _{i \in I} G_{i}$
might be isomorphic (in a natural sense) to $\sum _{i \in I} G_{i}$,
the game has no loop
(the relevant point in the above remark is that we consider
$\sum _{i \in I} G_{i}$  as a function
from some set $I$ to a set of combinatorial games;
were we considering $\sum _{i \in I} G_{i}$ as a set of options,
infinite runs would arise).

The above discussion introduces a natural notion of \emph{run}
on  $\sum _{i \in I} G_{i}$, but does not lead to a definite
notion of \emph{alternating run}. If we intend an alternating
run in the sense that First is always the first player to move
at each limit step, we get quite counterintuitive conclusions, so that 
we have dubbed such games as ``Dadaist''.
For example, if $I$ is countable and infinitely many
summands are $1-1$, then the resulting infinite game is always second 
winner, no matter the value of the remaining games.
So, for example, with the above notion of alternating run,
$ \omega + (1-1) + \omega ^ \omega + (1-1) + \omega ^ {\omega ^ \omega } +
(1-1) + \dots $ is a second winner game. 

The above fact can be proved as follows.
Let Second associate to each $G_i$ an infinite family of the 
$(1-1)$'s in such a way that these families are pairwise disjoint. 
If, whenever First plays on some $G_i$, Second replies on one of the associated games, Second never remains short of moves, since plays on
combinatorial games are finite. Hence First eventually loses, since she is always
the first at limit steps.

However, there are many possible different definitions
for an alternating run. Each  such definition
has some good characteristics, though we do not know
whether there is some optimal choice.
 
\emph{Regular alternating
runs.} A possibility is, at each limit step, to compute alternating runs by 
not considering pairs of consecutive moves 
by distinct players on the same
summand. Consider for example the sum 
$ 1 + (1-1)  + (1-1) +
(1-1) + \dots $, assume that
 Left first moves on the first $1$ and
subsequently she always replies
on the same summand where Right has just moved.
That is, after the first move,
 both players move alternately on each $1-1$
summand. Under the above rule, at each limit step,
the only remaining move to be considered is Left's first move.
 Then
Right is assumed to move first at the $ \omega$th step
and at subsequent limit steps, so that Left eventually wins.
With the above definition, Conway zero games,
that is, second winner summands, do not contribute
to the outcome of the sum and some aspects of the notion
are less counterintuitive. However, a countable sum 
which contains infinitely many positive numbers and
infinitely many negative numbers is always second winner,
no matter the values of the above numbers, and no matter the
forms of the remaining games. 

\emph{Games with a partition.}
 Another possibility is to partition $I$ into finitely many classes
and, at each limit step, consider some partition to be \emph{active}
if cofinally many moves have been done on games indexed by 
elements in the partition. At the limit step the player to move is the
first player who has ever moved \emph{on an active partition}.
With this definition we can recover parts of the theory of 
Conway equivalence. 
Indeed, with the above convention, let us 
``sum'' Dadaist games in the natural way
(by putting together all the summands) 
and consider those $G_i$s in some Dadaist game to belong
to a single partition.
For example, if 
$ G = 1/2$  
and 
$\mathbf H =   1-1  + 1-1 +
1-1 + \dots $,
then $G +\mathbf H =  1/2 + 1-1  + 1-1 +
1-1 + \dots $ with two partitions, the first partition is the 
singleton containing the first summand, and the other partition 
 contains the remaining games.
Under the above conventions, 
Left wins moving first on  $1/2$;
the second partition is the only one which remains active,
hence Right should move first at each limit step, so he loses.
Left wins moving second, as well.
If Right plays on $1/2$, Left can reply on the same game,
so Right is the first player to move on 
the partition which will turn out to be active.

So, with the above convention, $ 1-1  + 1-1 +
1-1 + \dots $ somewhat behaves as a Conway zero game, 
by taking sums (that is, fixing 
partitions) in the appropriate way.

\emph{Annihilating zeros.}
Finally, another possibility is the following.
At each limit step, consider only those moves
which have been played on summands which 
have not become $0$ at the limit step.  
Counting only the moves which remain,

  \begin{itemize}   
\item 
if both players have made cofinally infinite
moves, we require First to move at the limit step;

\item 
if, from some point on, all the moves are made by Some player,
we require the Other player to move at the limit step;

\item
if finitely many moves remain 
(or, more generally, the sequence of the remained moves
is indexed by a successor ordinal)
and Last 
has made the last move in the list, then
we require the Other player to move at
the limit step (this involves First moving,
if no move has remained). 
  \end{itemize}  

With the above rule, 
$ 1-1  + 1-1 +1-1 + \dots $
 behaves as a Conway zero game only in certain situations.
For example,
$1/2 + 1-1  + 1-1 +1-1 + \dots $
still remains a Second winner game.
If Right plays first on $1/2$,
it is enough that Left plays on the remaining
subposition $1$, turning the first component to $0$.  
All moves are then discarded at limit steps
and Right should move first at limit
steps, so he loses. If Left plays first,
then all moves are discarded, hence 
Left should move first at each limit
step, so she loses.

On the other hand,
$2 -1  + 1-1 +1-1 + \dots $
is Left winner.
The nontrivial case is when Left moves first.
Let her play on $2$, turning it to $1$.
She has always the possibility of moving
on the $1$ following the $-1$ on which
Right moves, so the first summand remains untouched.  
At limit steps, the only move which is  
 not discarded is Left's, hence, again,
Right should move first at each limit
steps, so he loses.

Note that the argument applies even to 
$(1 + \varepsilon ) -1  + 1-1 +1-1 + \dots $,
with $\varepsilon$ a positive infinitesimal number:
it is enough that Left has the 
possibility of moving on the summand 
without turning it to $0$ and in such a way that 
Right cannot move on the resulting option. 

As another example, Left 
always wins on 
$1-1/2+1-1/2+1-1/2+1- \dots $.
For example, if Left moves first, let her move on the first
$-1/2$ and subsequently on the first 
$-1/2$ available. If Right always replies on the same game,
then at the $ \omega$th step
$1+1+1+ \dots$ remains, so Left wins.  
Otherwise, some $-1$ (the remnant of  
some $-1/2$) persists, and Left's move on 
such a game is counted.
On the other hand, Right's moves always turn a game to $0$,
hence his moves are never counted.
In this situation Right must move first at
each limit step, and he loses in a game
with infinitely many $1$ summands   
and further $-1$ summands.

\smallskip 

Of course, we can also consider 
various combinations of the 
above possible rules
for defining an alternating run.

\smallskip 

Many problems remain; in particular,
for each of the above definitions of an alternating run,
 we do not even
know whether every such ``Dadaist'' game has a well-defined outcome,
namely, whether it is always the case that some player
has a winning strategy. Recall that we
are assuming the Axiom of Choice;
the situation in a choiceless setting appears to be very delicate
and  the very same definition of a winning strategy
requires some ingenuity.  

 We hope to be able to present further details in \cite{dada}.

\section{Appendix} \labbel{appen}

In the present appendix we describe an infinite sum
$\sumbul _{i  \in \mathbb N}$ 
which is Conway invariant, but has some undesirable 
properties.

\begin{definition} \labbel{sumpal} 
A play on 
$\sumbul _{i  \in \mathbb N} G_{i}$ goes as follows.
If First has no move on any $G_i$, First looses the game.
Otherwise, First chooses some natural number $n$
and she makes a move $G^F_i$ on some $G_i$  with $i \leq n$. 
The resulting position on $\sumbul _{i  \in \mathbb N} G_{i}$  
after  First move is 
 $G_0+ \dots + G^F_i+ \dots + G_n + G _{n+1} + \dots  $.

Then Second can either
  \begin{enumerate}[$\bullet$\phantom{$\bullet$}]  
 \item 
make a move on some $G_j$  with $j\leq n$
(of course, moving on $G^F_i$, not on $G_i$, if $j=i$), or
 \item[$\bullet\bullet$]
 choose some natural number $m \geq n$ and   
make a move $G_j^S$ on some $G_j$  with $j\leq m$
(with the same provisions as above).
  \end{enumerate} 

If Second chooses the eventuality $\bullet\bullet$,
the play continues on the finite sum
 $G_0+ \dots + G^F_i+ \dots + G_n + \dots +  G _{j}^S + \dots + G_m $
(possibly with the double subposition addendum $G^{FS}_i$, instead).
   
In the  eventuality $\bullet$, First can go on moving on games 
with index $\leq n$, but she cannot change the value of $n$.  
Then, again, Second has the possibility of  
keeping himself on indexes $\leq n$,
or of choosing some $m \geq n$
as in $\bullet\bullet$.
As soon as Second chooses 
$\bullet\bullet$, the play continues on a finite sum
$H_0 + \dots + H_m$, where  
 $H_0, \dots, H_m$ are the subpositions chosen so far
from the games $G_0, \dots, G_m$.

For short, each player has 
(exactly one time)
the possibility
of choosing some index, and then she 
has to move on some game with the same or  a smaller index.
After both players have made their choices
of the indexes, the play continues on the
finite sum of the  games
up to the larger chosen index.
The same description applies to \emph{runs},
ie, when the same player is supposed to make
 two or more consecutive moves.
Notice that, after the choice of $n$ by First
and until Second applies the eventuality
$\bullet\bullet$, the play essentially proceeds on
$G_0+ \dots +G_n$, which is  a finite sum of games,
hence either the game ends after a finite number
of moves, or Second is forced to apply
$\bullet\bullet$. After this, the game 
definitely proceeds 
on a  finite sum with a fixed index $m$, hence 
the game eventually terminates.
Thus 
$\sumbul _{i  \in \mathbb N} G_{i}$
has no infinite run.

We now define 
$\sumbul _{i  \in \mathbb N} G_{i}$  formally.
We need to introduce auxiliary games
such as  
$ \left( \sumbul _{i  \in \mathbb N} H_{i}\right ) / {L,n}$,
with the intended meaning that ``Left has
already made the choice of her index 
and she has chosen $n$''. 
For clarity, we may possibly write
 $(H_0+H_1+H_2+ \dots )^\bullet / {L,n}$ in place of
$ \left( \sumbul _{i  \in \mathbb N} H_{i}\right ) / {L,n}$.
Thus
\begin{align*}  
&\big( \sumbul _{i  \in \mathbb N} H_{i}\big) / {R,n} 
=
\\
&\big\{ \, (H_0+ \dots +H^L_i+ \dots + H_n + \dots )^\bullet / {R,n}, \ 
H_0+ \dots +H^L_j+ \dots + H_m  \  \big| 
\\
&  (H_0+ \dots +H^R_i+ \dots + H_n + \dots ) ^\bullet / {R,n}  \, \big\}, \text{ and}
\\
&\big( \sumbul _{i  \in \mathbb N} H_{i}\big) ^\bullet/ {L,n} 
=
\big\{ \, (H_0+ \dots +H^L_i+ \dots + H_n + \dots ) ^\bullet/ {L,n} \  \big| 
\\
&  (H_0+ \dots +H^R_i+ \dots + H_n + \dots ) ^\bullet/ {L,n}, 
\ 
H_0+ \dots +H^R_j+ \dots + H_m \ \big\}
\end{align*} 
where $m$ is always intended to be $\geq n$.  
As in \ref{sum}, the  definition 
is by transfinite induction on the natural sum of the 
birthdays of $H_0, \dots, H_n$.

Finally,
\begin{align*}
 \sumbul _{i  \in \mathbb N} G_{i} =
& \{ \, (G_0+ \dots +G^L_i+ \dots + G_n + \dots ) ^\bullet / {L,n}\  \big | \ 
 \, 
\\ & 
(G_0+ \dots +G^R_i+ \dots + G_n + \dots )^\bullet / {R,n} \} 
  \end{align*}     
\end{definition}   

A possible variation on the above definition requires
that, until Second applies 
$\bullet\bullet$, she must move on the same game
First has just moved. 
As in Theorem \ref{equiv} below,
Conway invariance is maintained
with this variation; anyway   paradoxical properties remain.

\begin{observation} \labbel{obspal}   
Suppose that $(G_i) _{i  \in \mathbb N} $
is a sequence of games.
  \begin{enumerate}    
\item
If each $G_i$ is impartial (a number, dicotic, $>0$, $\geq 0$,  $<0$, $\leq 0$),
then  $\sumbul _{i  \in \mathbb N} G_{i}$ is
impartial (a number, dicotic, $>0$, $\geq 0$,  $<0$, $\leq 0$).

\item
If there are infinitely many partial sums which are
$\geq 0$ (resp., $\leq 0$), then Left (resp., Right)
wins when playing second on  $\sumbul _{i  \in \mathbb N} G_{i}$.
 In particular, if both eventualities hold, then 
$\sumbul _{i  \in \mathbb N} G_{i}$ is second winner.

\item
Similar to the case of  $\sum _{i  \in \mathbb N} G_{i}$,
 if  each $G_i$ is an ordinal, then
 $\sumbul _{i  \in \mathbb N} G_{i}$ is the supremum
(in the ordinal sense)
of the partial sums.

\item
 An infinite sum of 
  impartial games is a 
second winner game, if infinitely many summands
are first winner.

If only finitely many summands
are first winner, the outcome is given by the finite sum
of such first winner summands, by Corollary \ref{cor}
to be proved below,
since an impartial game is either first  or second winner.
  \end{enumerate}
\end{observation} 

\begin{proof}
(1) - (3) are proved as in Observation \ref{obs}. 

(4) Let First play on such a sum and choose some $n$.
Second can always move on $G_0+ \dots + G_n$
until all such games are  $0$, since the games are impartial.
If Second has won, we are done.

Otherwise, just let him chose the first $m>n$
such that $G_m$ is first winner (by assumption, there are infinitely
many such games). Now Second has a winning
move on the first winner game
$G_m$, hence he wins, since the remaining
games are second winner and $G_0+ \dots + G_n$
has been emptied. 
\end{proof}

\subsection{Invariance under Conway equivalence} \labbel{inva}

It is standard to see that
$- \sumbul _{i  \in \mathbb N} G_{i}$ 
is
$\sumbul _{i  \in \mathbb N} (- G_{i})$.

\begin{theorem} \labbel{equiv}
If $G_i$ is Conway equivalent to
$H_i$, for every $i \in \mathbb N$,
then   $\sumbul _{i  \in \mathbb N} G_{i}$ 
is Conway equivalent to
$\sumbul _{i  \in \mathbb N} H_{i}$. 
 \end{theorem} 

\begin{proof}
In view of the comment before the statement,
we have to show that 
$\sumbul _{i  \in \mathbb N} G_{i} + 
\sumbul _{i  \in \mathbb N} (- H_{i})$
is second winner. 

Suppose that, say, Right moves first 
and decides to make his first move on the sum
on the left (the other cases are treated symmetrically).
So Right chooses some $n$
for the game $\sumbul _{i  \in \mathbb N} G_{i}$ 
together with  some option $G_i^R$. 
Since, by assumption, $G_i$
is Conway equivalent to $H_i$,
$G_i - H_i$ is second winner, hence
Left has a winning move on   
$G_i^R - H_i$.
There are two cases.

(a) If the winning strategy for Left on $G_i^R - H_i$
prescribes a move on 
$-H_i$, let Left choose the same $n$ as above
on   $\sumbul _{i  \in \mathbb N} (- H_{i})$
and make the prescribed move on $-H_i$.
At this point, Right cannot change the $n$
on the left-hand sum, but he can eventually change
the $n$  on the right-hand sum.
Until he makes this latter change, the play continues on
(the remaining options of)
$G_0 + \dots + G_n - H_0 - \dots H_n$,
on which Left has a winning strategy. 
As far as Right changes $n$ to some  $m$
on the    right-hand sum, consider that
  Left has a winning strategy on
$G_{n+1} + \dots + G_m - H_{n+1} - \dots H_m$.
It is thus  enough that Left chooses $m$
for $\sumbul _{i  \in \mathbb N} G_{i}$ 
as soon as she moves on 
$G_0 + \dots + G_m$,
in order to get a winning strategy  for
$\sumbul _{i  \in \mathbb N} G_{i} + 
\sumbul _{i  \in \mathbb N} (- H_{i})$.

(b) Otherwise, after Right's first move,
the winning strategy for Left prescribes 
a move on $G_i^R$. In this case Left
is not required to choose some $m \geq n$,
according to the possibility given by $\bullet$
in Definition \ref{sumpal}, so 
she
simply  moves on $G_i^R$. 
As long as the players go on moving 
on the $G_j$s, Left continues with 
her winning strategy on $G_0 + \dots + G_n - H_0  \dots - H_n$,
 never changing the value of $n$. 
Sooner or later, some player will
be forced to  move on the $H_j$s.

(b1) 
If Left is the first to move on the $H_j$s, 
she can win by applying the same strategy as described in (a), namely
 choosing $n$ in  $\sumbul _{i  \in \mathbb N} (- H_{i})$.

(b2) 
If Right is the first to move on the $H_j$s, 
he must choose some $n'$ for $\sumbul _{i  \in \mathbb N} (- H_{i})$. 
If $n' \leq n$, let Left choose 
$n$ on  $\sumbul _{i  \in \mathbb N} (- H_{i})$
as soon as she moves there.
If $n' > n$, let Left choose 
$n'$ on  $\sumbul _{i  \in \mathbb N} (G_{i})$
as soon as she moves there.
In any case, the game continues on 
(the remaining options of)
$G_0 + \dots + G_m - H_0  \dots - H_m$,
for $m= \max \{ n,n' \} $ and, by the assumptions,
Left has a winning strategy on this game. 

We have described a winning 
strategy for the second player on 
$\sumbul _{i  \in \mathbb N} G_{i} + 
\sumbul _{i  \in \mathbb N} (- H_{i})$
and this means that 
 $\sumbul _{i  \in \mathbb N} G_{i}$ 
and
$\sumbul _{i  \in \mathbb N} H_{i}$
are Conway equivalent. 
\end{proof}  

\begin{remark} \labbel{serve}
Allowing to Second the possibility of using the rule
$\bullet$ in Definition \ref{sum}  is necessary for Theorem \ref{equiv}  to hold.
Namely, Theorem \ref{equiv} does not generally hold if,
as in Definition \ref{sum},
Second is required to choose some $m \geq n$,
as soon as he moves on   $\sum _{i  \in \mathbb N} G_{i}$.

For example, $0+1+1+1+\dots $ is Conway equivalent
to $ \omega$; in fact, Right cannot even move
 on  $0+1+1+1+\dots $, while for Left a move
on  $0+1+1+1+\dots $ corresponds to choosing some 
$n-1$, as a game. 

On the other hand, if Second is not allowed to move according to
$\bullet$, then $(-1+1)+1+1+1+\dots $
is strictly smaller than $ \omega$.
Indeed, Right has a winning strategy on the sum
of $(-1+1)+1+1+1+\dots $ and $- \omega$.
Just let Right move on $(-1 +1)$ (no matter the $n$
he chooses), so Left should choose some $m$   
on $1+1+1+1+\dots $, which is losing against
$- \omega$.
In fact, without allowing $\bullet$, $(-1+1)+1+1+1+\dots $
is Conway equivalent to $ \omega - 1$.
 \end{remark}   

\begin{corollary} \labbel{cor}
Second winner games do not influence a sum 
 $\sumbul _{i  \in \mathbb N} G_{i}$.
This means the following.

If all but finitely many $G_i$
are second winner,
then  
 $\sumbul _{i  \in \mathbb N} G_{i}$ and
$ G_{i_1} + G_{i_2} + \dots + G_{i_h} $
are Conway equivalent, where 
$ G_{i_1} $, $   G_{i_2} $, \dots, $  G_{i_h} $
lists those games which are not second winner.

If there are infinitely many $G_i$
which are not second winner, let $( H_n) _{n \in \mathbb N} $  
 be the  reindexing (respecting the order)
of the set of those $G_i$ which are not second winner. Then  
 $\sumbul _{i  \in \mathbb N} G_{i}$ and
$\sumbul _{n  \in \mathbb N} H_{n}$
are Conway equivalent.
\end{corollary}

 \begin{proof} 
By Theorem \ref{equiv}, 
$\sumbul _{i  \in \mathbb N} G_{i}$ is Conway equivalent
to the series in which all second winner games are changed
to $0$. Then the corollary is an immediate application of  
a mirror-image strategy.
\end{proof}

As in Remark \ref{nolim}(b),
we can similarly define a permutation invariant sum, by letting
each player choose a subset of $\mathbb N$.
The same proof as of Theorem \ref{equiv}
shows that this sum is Conway invariant.  

Similarly, we can transfinitely iterate the sum defined in 
Definition \ref{sumpal} as in 
Remark \ref{nolim}(f).
In more detail, for $\alpha$
a limit ordinal, First chooses some $\beta < \alpha $ and moves on
some $G_ \gamma $ with $\gamma < \beta $.
Second can either play on some 
$G_ {\gamma'} $ with $\gamma' < \beta $,
or (sooner or later) choose some $\delta$ such that  
$\beta \leq \delta < \alpha $.
After Second makes the latter choice,
 the play definitely continues on (the remaining
options of) $\sum ^\bullet _{ \beta  < \delta  } G_{ \beta }$.

A result analogue to Theorem \ref{equiv}
holds for the above transfinite $\sum ^\bullet_{ \beta  < \alpha } G_{ \beta }$.

\subsection{Paradoxical sums} \labbel{parad} 

The main novelty in the present approach seems
to be the possibility of considering sums of games which are not 
numbers. However, Observation
\ref{obspal}(4) shows that  there are situations in which 
the outcome of the sum is trivial.
Moreover, there are quite paradoxical results, as shown
by the next proposition.

\begin{proposition} \labbel{small}
Suppose that $(G_i) _{i  \in \mathbb N} $
is a sequence of games.

  \begin{enumerate}   
 \item   
If each $G_i$ is  $>0$, then  
  $\sumbul _{i  \in \mathbb N} G_{i} \geq 1$.
 
\item 
If each partial sum is $<1$,
then 
$\sumbul _{i  \in \mathbb N} G_{i} \leq 1$.
\item
Thus, if both assumptions in (1) and (2) apply,
then 
$\sumbul _{i  \in \mathbb N} G_{i} = 1$.
\end{enumerate} 
\end{proposition} 

 \begin{proof}
 (1) We have to prove that Left always wins
when playing
second on the sum of  $\sumbul _{i  \in \mathbb N} G_{i}$ and $-1$.

If Right moves first, he loses moving on $-1$,
since all the $G_i$ are  $>0$. 
So let Right move on the infinite  sum, choosing some index $n$.
Since all $G_i$ are $> 0$, Left wins on
$G_0+G_1+ \dots + G_i^R+ \dots+ G_n$.
 Then Right cannot change his index $n$,
so he should move on $-1$, turning it to $0$.  
Since the $G_i$ are $>0$,
Left wins by choosing any $m > n$.

(2) We have to prove that Right always wins
when playing
second on the sum of  $\sumbul _{i  \in \mathbb N} G_{i}$ and $-1$.
Indeed, if Left moves first, she should move on the infinite sum,
choosing some index $n$.
By the assumption, 
Right wins on
$-1 + G_0+G_1+ \dots + G_i^L+ \dots+ G_n$.
Left cannot change her index $n$ and
Right is not forced to choose a larger index,
so Right wins on the main game.
Of course, this proof of (2)
applies also to $\sum _{i  \in \mathbb N} G_{i}$,
while (1) does not hold for $\sum$, since in the case of
$\sum $ Left is forced to choose some index $m$
as soon as she plays on the sum.  
 \end{proof} 

By Proposition \ref{small},
we get the unwanted result that, 
according to Definition \ref{sumpal},
$1/4+1/8+1/16+ \dots =1$.
Even worse, a sum of infinitely many 
positive infinitesimal games is still $1$
(the arguments in the proof of Theorem \ref{reals} 
do not work for $\sumbul$, since in step (b1)
in the proof Right is forced to immediately
choose some $m$, while Right is allowed to  wait,
under the rules in Definition \ref{sumpal}).  

We have seen that the operation $\sumbul$ 
is so weird that it cannot be
 seriously considered as a natural operation.
However, possibly, the ideas in Definition \ref{sumpal}
might be modified in order to obtain some Conway
invariant operation with less paradoxical
properties.

\end{document}